\documentclass[12pt]{amsart}
\usepackage{geometry}                
\usepackage[parfill]{parskip}    
\usepackage{graphicx}
\usepackage{amssymb}
\usepackage{epstopdf}
\usepackage{color}

\newtheorem{theorem}{Theorem}[section]
\newtheorem{prop}[theorem]{Proposition}
\newtheorem{lemma}[theorem]{Lemma}
\newtheorem{cor}[theorem]{Corollary}

\newcommand{\btheorem}{\begin{theorem}}
\newcommand{\etheorem}{\end{theorem}}
\newcommand{\bprop}{\begin{prop}}
\newcommand{\eprop}{\end{prop}}
\newcommand{\blemma}{\begin{lemma}}
\newcommand{\elemma}{\end{lemma}}
\newcommand{\bcor}{\begin{cor}}
\newcommand{\ecor}{\end{cor}}

\newcommand{\bg}{{\mathfrak g }}
\newcommand{\ba}{{\mathfrak a }}
\newcommand{\bn}{{\mathfrak n }}

\newcommand{\bk}{{\mathfrak k }}
\newcommand{\bq}{{\mathfrak q }}

\newcommand{\bh}{{\mathfrak h }}
\newcommand{\bt}{{\mathfrak t }}
\newcommand{\bff}{{\mathfrak f}}

\newcommand{\bp}{{\mathbf p}}

\newcommand{\bR}{{\mathbb R}}

\newcommand{\bC}{{\mathbb C}}
\newcommand{\bZ}{{\mathbb Z}}

\newcommand{\R}{{\mathbb R}}

\title{A Plancherel formula for $L^2(G/H )$ for almost symmetric subgroups}
\author{Bent {\O}rsted and Birgit Speh}
\address{B.{\O}rsted, Department of mathematics, Aarhus University, 8000 Aarhus C, Denmark
\newline
Email:
orsted@math.au.dk}
\address{B. Speh, Department of Mathematics, Cornell University, Ithaca 
  NY 14853, USA.  
\newline
Email: 
speh@math.cornell.edu}
\thanks{Research by B.Speh partially supported by NSF grant DMS-0901024 }

\begin{document}

\maketitle

\begin{abstract}
In this paper we study the Plancherel formula for a new class of homogeneous
spaces for real reductive Lie groups; these spaces are fibered over non-Riemannian
symmetric spaces, and they exhibit a phenomenon of uniform infinite multiplicities.
They also provide examples of non-tempered representations of the group
appearing in the Plancherel formula. Several classes of examples are given.

\end{abstract}



\section{Introduction:}
Considerable efforts have been devoted to obtaining the Plancherel formula for
homogeneous spaces of the form $G/H$ with $G$ a real reductive Lie group and
$H$ a symmetric subgroup, a program completed by T. Oshima, P. Delorme,
E. van den Ban, and H. Schlichtkrull. This is a central theme in harmonic analysis,
and there are a number of natural ways to extend such a program. One is
to consider spherical spaces, i.e. where the homogeneous space admits an open
orbit of a parabolic subgroup. In this paper we shall rather extend the interest
to \begin{enumerate}
\item square-integrable sections of homogeneous line bundles over symmetric
spaces, and to \item spaces fibered over symmetric spaces. 
\end{enumerate}

Of course, these two
questions are related, and we shall find several classes of spaces where rather
explicit answers can be found. As an example consider $G = SL(2,\R)$ with
$H$ the connected diagonal subgroup; for each unitary character of $H$
we may consider the space (1) and the corresponding Plancherel formula:
This turns out to be independent of the character, and hence the space as in (2) above
(in our case just the group $G$) has the same $L^2$-content as the symmetric
space, only with infinite multiplicity. It is perhaps a little surprising, that one may
thus find embeddings of  \linebreak e. g. the discrete series of $G$ in a uniform way in each
of the spaces of sections (1).

\medskip
To be more specific, our motivation  in undertaking this work was to understand  the disintegration of the representation of  a semi-simple Lie group $G$ on the space $L^2(G/H_{ss})$  where $H_{ss}$ is a semi-simple subgroup which differs from  a symmetric subgroup by a non-compact central real factor. 
In this paper we  study this situation for the simplest non-symmetric subgroups $H_{ss}$  from the point of view of harmonic analysis and obtain a Plancherel theorem for space $L^2(G/H_{ss})$ in terms of the one for $G/H$.
 
 \medskip
Recently Y. Benoist and T. Kobayashi proved general criteria to determine if for a  semi-simple subgroup $H$  the spectrum of $L^2(G/H)$ contains non-tempered representations; this they use to determine in many examples if $L^2(G/H)$ is tempered. Here a representation is called  tempered if it appears in the usual Plancherel formula for $L^2(G)$. However
 these authors do not obtain any results concerning the multiplicity of the representations in the Plancherel formula. 
 By obtaining a Plancherel formula for $L^2(G/H_{ss})$ we are in a position to determine exactly in our examples which non-tempered representations appear in the spectrum, and also to show that they appear with infinite multiplicities. 
 
 \medskip
 
We consider a noncompact subgroup $H = H_{ss}Z_H$ where $H$ is a subgroup of finite index in the fixpoints of an involution of $G$  and $Z_H \simeq \bR$ is a subgroup of finite index of the center of $H$. Under these assumptions we show that 

\bigskip 
\noindent 
 {\bf Theorem} \ {\it  As a left regular  representation of $G$}
\[ L^2(G/H_{ss} )   \simeq   L^2(G/H) \otimes L^2(Z_H).\]

\bigskip 
It is instructive to compare with the situation where the central subgroup is compact, e.g. the case
of $G$ a simple non-compact Lie group and $K$ a maximal compact subgroup with a one-dimensional center
$Z$. Here $G/K$ is a non-compact Riemannian symmetric space of Hermitian type, and $L^2(G/K)$ has
a different Plancherel decomposition than $L^2(G/K, \chi)$, the square-integrable sections of the
line bundle induced from 
a non-trivial unitary character $\chi$ of $Z$. In particular the first space contains
no discrete series representations, whereas the second space typically does.
Compare with Proposition III.3 for our situation of a non-compact center.      

\medskip 
The paper is organized as follows:  In section II, we show that we can regard $H$ as a subgroup of finite index in the Levi subgroup of a parabolic subgroup with abelian nilradical. In the following section we prove our main theorem
above. In the last section we discuss some examples. In particular we note that we find several examples of
non-tempered homogeneous spaces, some of them new; quite possibly our method could extend to other instances of
Plancherel theorems, such as cases of vector bundles (as opposed to the cases of line bundles treated here).

\bigskip
The second author would like to thank the department of mathematics  of the University of Aarhus for their  hospitality.

\bigskip

\section{Notation and Preliminaries}

{\it In this section we  introduce the notation and prove some preliminary results;}

\subsection{Notation and assumptions}
Let $G$ be a  real linear semi-simple connected algebraic group with maximal compact subgroup $K$ and complexification $ G_\bC$. We consider $G,\mbox{ and } K \subset G_\bC$ as subgroups.

\medskip
\begin{prop} Suppose that $P=LN$ is a maximal parabolic subgroup with an abelian nilradical $N$.
Then $L$  is the fixed point set of an involution \[\tau: G \rightarrow G. \]
 \end{prop}
\proof (due to Dan Barbasch) We consider a maximal split Cartan subgroup and its corresponding 
complex Cartan subalgebra. A parabolic subalgebra is given by removing some simple roots from the
diagram. The only way to get an abelian nilradical is to remove a single
simple root which appears with coefficient at most one 1 if we write the roots as linear combinations of simple roots. 
The involution $\tau$ is then  conjugation by  $ exp({i\pi\varpi})$ where  $\varpi$ is the coroot of
the simple root which was removed.      \qed

\medskip
Let $H$ be a subgroup of the Levi subgroup  $L$ of $P$ which contains the connected component $L^0$ of $L$. 
Then $H=H_{ss}Z_H$ where $Z_H$ is a one dimensional connected subgroup in the center of $H$  and $H_{ss}$ is semisimple or discrete.
 
\medskip 

\noindent
{\bf Example 1: }   $G=SL(2,\bR)$, $L$ diagonal matrices which are the fixpoints under the conjugation by the diagonal matrix of order 2 and determinant -1. Alternatively we consider the adjoint representation. Then $L$ is the stabilizer of a semi-simple nontrivial element of order 2.  It is also the fixpoint set of the automorphism by the adjoint acton of the matrix
  \[
\left(
\begin{array}{cc}
i  & 0    \\
0  & -i    
 
\end{array}
\right)
= exp (\pi i \left( \begin{array}{cc}
1/2 & 0 \\
0 & -1/2 
\end{array}
 \right) \]
We have to consider  2 subgroups $H$ and $H_{ss}$
   \newline a.) $H=L$,  $H_{ss }= \bZ_2$, $Z_H  =\bR_+$ . $G/H_{ss} =PSL(2,\bR)$ .\newline
  b.) $H=L^0$, $H_{ss}= I$   and  $G/H_{ss} =SL(2,\bR)$.

\medskip

\begin{prop} 
Suppose that $F$ is the fixpoint set of an involution $\tau: G \rightarrow G $. Assume in addition that it is a product  $F=F_{ss} Z_F$ where $Z_F$ is a subgroup of the center of $F$ isomorphic to $\R^+$ and $F_{ss} $ is a semi-simple group. Then $F$ is contained in the Levi subgroup of a maximal parabolic subgroup P with abelian nilradical $N_P$.
\end{prop}
\proof We choose maximally  split Cartan subgroup $C \subset F$ with complexified Lie algebra $\bh_\bC$. We choose the simple roots of $\bh_\bC, \bg_\bC$ so that they are simple roots in $\bh_\bC, \bff_\bC$. ($\bff$ before $\bg$,
$\bff$ the Lie algeba of $F$). Then $\bff$ is the Levi subalgebra  of a maximal parabolic subalgebra $\bp_\bC = \bff_\bC \oplus \bn_\bC$ of $\bg_\bC$.  

It remains to show that N is abelian.
Since $\tau$ leaves $C$ invariant the induced homorphism of $\tau:N \rightarrow N$  is equal to -1.  Since $\tau$ induces a Lie algebra homorphism and hence preserves the Lie bracket in $\bn_\bC$, the results follows from the observation that $ \tau (X) = -X $ and $\tau (Y) = -Y$ then  $\tau ([X,Y]) = [X,Y]$.
\qed

Note that in the setting above, we have a direct product
decomposition $G/H_{ss} \times Z_H = G/H $. This will be useful later in connection with integration
over this space, and in considering the corresponding $L^2$ - space.

\medskip 
\noindent

\subsection{About $L^2(G/H)$} 

\noindent
Keep the assumptions on $G, H, Z_H$ as above.
We extend a unitary character $\chi \in \widehat{Z_H} $ to a character of $H$ and consider 
the unitary induced representation $\mbox{Ind}_H^G \chi$ on $L^2(G/H)_{\chi ^{-1}}.$ Normalize Plancherel measures
on $Z_H$ and its dual group in the usual way.

\begin{prop}
As a representation of $G$ 
\[ L^2(G/H_{ss}) = \int_{\chi \in \widehat{Z_H}}L^2(G/H)_{\chi ^{-1}} d\chi\]
\end{prop}
 \proof
For $f \in L^2(G/H_{ss})$ and $\chi \in \widehat{Z_H} $ define 
\[F(\chi, g) = \int _{Z_H} f(gz)\chi(z)^{-1}dz \]
Then for $z_0 \in Z_H$ \[F(\chi, gz_0) = F(\chi,g)\chi^{-1}(z_0) \] 
so $F(\chi) \in L^2(G/H)_{\chi ^{-1}}$.
By Fourier analysis on $Z_H$ we have 
\[\int _{\chi \in \widehat{Z_H} } |F(\chi,g) |^2 d\chi = \int_{z\in Z_H}|f(gz)|^2 dz.\]
So
\[ \int_{G/H_{ss} }|f(g')|^2dg' =\int_{G/H} \int _{Z_H}|f(g'z)|^2dzdg'\]
completes the proof.   \qed

\section{Main Results}

\medskip
{ \it In this section we  relate the Plancherel formula for  the left regular representation of $G$ on $L^2(G/H_{ss} )$ to the Plancherel formula for the left regular representation on $L^2(G/H)$. It turns out that these two
spaces have the same content of unitary representations of $G$, only differing by their multiplicities.}

\subsection{Induction to the parabolic subgroup P}

\begin{lemma}
\noindent
Let $\widehat{N}$ the dual group of $N$. There exist finitely many open $H$ orbits $\mathcal O_i$ 
in $\widehat{N}$ so that 
$ \widehat{N}$ is the closure of their union $\bigcup_i { \mathcal O_i}$. 
\end{lemma}
\proof 
Here we refer to results by Wallach, see
[W].  
Here he proves that our parabolic algebras are "very nice" since they have abelian nilradicals
- see in particular Corollary 6.4. In particular there is only one open orbit of $L$ on $N$.\\
Since our group $H$ is a subgroup of finite index in $L$ we will get a finite number of open orbits
with dense union. Actually, the statement that "open orbit is generic" (i.e. "nice
parabolic") would suffice for
our purposes here.
\qed


\bigskip
Let $\chi \in \widehat{Z_H}$. We consider again $\chi$ as a character of $H$ and consider again
the unitary induced representation
$\mbox{Ind}_H^P \chi $.

\begin{prop} 
Let $\chi$ and $\tilde{\chi}$ be  unitary characters of $Z_H$ considered as characters of $H$. Then 
(equivalence of representations)  
\[Ind_H^P \chi = Ind_H^P \tilde{\chi}.\]
\end{prop}
\proof  
We denote the induced representations
act on functions $F \in L^2(N)$ by
    \[ \rho_\chi (n_0)F(n) = F(n\cdot n_0)\]
    \[ \rho_\chi (h_0) F(n) = \chi(h_0)F(h_0^{-1} n h_0)\]
Using the Fourier transform we realize the representation $Ind_H^P \chi$ on $L^2(\hat{N})$ . It is a direct sum of irreducible representations on $L^({\mathcal O}_i)$ where
 \[\hat{ \rho}_\chi (n_0) \mbox{ is a multiplication operator }\]
 \[\hat{ \rho}_\chi (h_0){\hat F}(\xi) = \chi(h_0)J(h_0^t \xi)^{1/2} \hat{F}(h_0^t\xi) \]
 The other representation on the orbit is obtained by multiplication of the right hand side of the second equation with a character $\chi_1=\tilde{\chi}\chi^{-1} $ of $H$.
 In each orbit we fix and element $\xi_i$.
 
 We get a intertwining  operator on each of the irreducible representations  by
    \[ I( F )(\xi) = \chi_1(\xi) F(\xi)\]
Here $\xi = h\xi_i$ and $\chi_1(\xi ) :=\chi_1(h) $.

\qed

\medskip
\noindent
{\bf Example 2:} Consider the group $P=HN$ with 
  \[ H= \{
\left(
\begin{array}{cc}
a & 0    \\
0  & 1   
 
\end{array}
\right),\ \ |  \  a>0 \}\]
and
  \[ N= \{
\left(
\begin{array}{cc}
1  & b    \\
0  & 1   
 
\end{array}
\right)\ \ | \ b \in \bR \} \]
We note that there are 3 orbits of H on
 \[
\hat N = \{\xi_t \   | \ \xi_t(\left(
\begin{array}{cc}
1  & b    \\
0  & 1   
 
\end{array}
\right)  )= e^{it \cdot b}\}.
\]
namely ${\mathcal O}^+= \{\xi_t \ |\ t>0 \}$ , ${\mathcal O}^-= \{\xi_t \ |\ t<0 \}$ and ${\mathcal O}^1= \{\xi_0 \}.$
The unitary representation  $\rho_1$ of $P$ induced from the trivial representation of $H$  acts on $L^2( N)$ by
\[ \rho_1(\left(
\begin{array}{cc}
a & 0    \\
0  & 1   
 
\end{array}
\right))F(x) = a^{1/2}F(ax) \]
and 
\[\rho_1(\{
\left(
\begin{array}{cc}
1  & b    \\
0  & 1   
 
\end{array}
\right))F(x) =F(x+b). \]
To analyze this representation we consider the Fourier transform of  $L^2( N)$. 
The representation is a direct sum of 2 unitary representations of functions whose Fourier transform has support in  $\xi \in  {\mathcal O}^+$ and in 
$\xi \in {\mathcal O}^-$

We consider  $\chi_s: a\rightarrow a^{is} $ as a character of $H$.    After applying the Fourier transform the  representation $\widehat{\rho_{s}}$  induced from $\chi_s$ has the form

\[\widehat{ \rho_s}(\left(
\begin{array}{cc}
a & 0    \\
0  & 1   
 
\end{array}
\right))\hat{F}(\xi) = a^{-1/2}a^{is} \hat{F}(a^{-1}\xi) \]
and 
\[\widehat{\rho_t}(
\left(
\begin{array}{cc}
1  & b    \\
0  & 1   
 
\end{array}
\right))\hat{F}(\xi) =e^{ib\xi}{\hat{F}(\xi)} \]
 The equivalence of the representations $\rho_s$ and $\rho_1$ follows from  the intertwining operator
\[ {\mathcal I}_s : \rho_0 \rightarrow \rho_s \ \ \ \ \mbox{defined by} \ \ \ \ \ \ 
{\mathcal I}_s\hat{F} (\xi ) = \xi^{is} \hat{F}(\xi ) . \]

\subsection{Induction to G}
\begin{prop} 
Let $\chi$ and $\tilde{\chi}$ be  characters of $Z_H$ considered as characters of $H$. As
representations of $G$ we have (equivalence)  \[Ind_H^G \chi = Ind_H^G \tilde{\chi}.\]
\end{prop}
\proof By induction by stages proposition III.2
 \[Ind_H^G \chi = Ind_P^G Ind_H^P \chi = Ind_P^G Ind_H^P \tilde{ \chi} = Ind_H^G\tilde{ \chi }\]
\qed

\medskip

\noindent
{\bf Example 3}  $G=SU(1,1)$ and \[H = A =\mbox{exp} 
\,\bR \left(
\begin{array}{cc}
  0& 1     \\
  1   & 0  
  
\end{array}
\right)
.\]
We identify $G/K$ with the complex unit disk ${\mathcal D}$. We realize the discrete series representations $D_n$ in the holomorphic functions on ${\mathcal D}$.
Then we have the $H$ - invariant distribution vector, giving the imbedding into $L^2(G/H)$, 
\[v^*= (1+z^2)^{-n/2}  \in D_n^{-\infty, H}\]
and similarly the distribution vector   
\[v^*= (1+z^2)^{-n/2}  (\frac{1-z}{1+z})^{i\lambda} \in  D_n^{-\infty, H ,\chi_\lambda}\]
transforming by the character $\chi_\lambda$ of $H$. So indeed every discrete series representation 
occurs in every $L^2(G/H)_{\chi_\lambda}$

\bigskip
\begin{theorem}  As a left regular  representation of $G$
\[ L^2(G/H_{ss} )   \simeq (Ind _H^G 1 )\otimes L^2(Z_H) \simeq  L^2(G/H) \otimes L^2(Z_H).\]
\end{theorem}
\proof This follows from propositions II.3 and III.3. \qed

\medskip 
\begin{cor} 
 All irreducible representations in the discrete spectrum of $L^2(G/H_{ss})$ have infinite multiplicity.
 \end{cor}

\noindent
{\bf Definition}
Following Benoist and Kobayashi we say that  $L^2(G/H_{ss})$ is not tempered if the representations in the Plancherel formula for the right regular representation of $G$ on $L^2(G/H_{ss}$ are not a subset of the representations of the Plancherel formula for $G$.

  \medskip 
\begin{cor} 
 $L^2(G/H_{ss})$ is tempered iff $L^2(G/H)$ is tempered.
 \end{cor}

\medskip
{\bf Example 1 (continued)} \ \   $G=SL(2,\bR)$, $H$ diagonal matrices, Then $X= G/H$ is a hyperboid and 
\[L^2(G/H) = \oplus_{\nu \mbox{ \it even }} D_\nu + 2 \int_0^\infty \pi_{it} \]
where $D_\nu $ are the discrete series representations with parameter $\nu$ and $\pi_{it}$ are the 
tempered spherical principal series representations with parameter $it$.  
Here $H_{ss }= \bZ_2$, then  $L^2(G/H_{ss}) =L^2(PSL(2, \bR ))$ and so the left regular representation contains the even discrete series representations with $\infty $ multiplicity.

If $H $ is connected, then $L^2(G/H)$ contains all discrete series representations and so does the left regular representation of G on $L^2(G)$.


 \bigskip
 \section{More Examples}
 {\it We discuss in this section some interesting examples of groups $G$ and $H_{ss}$,
illustrating our results; one aspect is to find reductive spaces that are not tempered.}

     We use the Plancherel formula  to determine if $L^2(G/H_{ss})$ is tempered. Some of our examples are also contained  in [B-K], where they are obtained with a different technique;  others are new. 
    
     E. van den Ban and H. Schlichtkrull proved a Plancherel Formula for $L^2(G/H)$ for a fixed point set $H$ of an involution $\tau $ of $G$.  They showed that only discrete  series representations of $L^2(G/H) $ and principal series representations unitarily induced from a $\theta \tau$ invariant parabolic  $MAN$, a discrete series representation $\pi$ of $M/M\cap H$ and a unitary character of $A$ contribute to the Plancherel formula. [B-S]
On the other hand    the work 
    of M. Flensted-Jensen and Oshima-Matsuki shows that the discrete spectrum of $G/H$ is nontrivial iff 
     \[ \mbox{rank } G/H  = \mbox{rank } K/K\cap H.\]
 and a parametrization of the representations in the discrete spectrum was obtained by T.Matsuki and T.Oshima [M-O]. See also [S]. We will make extensive use of these results in the proofs of our examples.

\noindent {\bf Remark 1 :}  Induction by stages enlarges the set of pairs $G, \tilde H$ so 
that $L^2(G/\tilde{H})$ is tempered. (See theorem 4.2 in [F]; here the point is that induction
preserves weak containment, so if we have groups $H \subset \tilde{H} \subset G$ so that
$L^2(G/H)$ is tempered and we know that $L^2(\tilde{H}/H)$ contains the trivial representation
weakly, then also $L^2(G/\tilde{H})$ is tempered.)

 \noindent{ \bf Remark 2 :} The non tempered representations in the discrete spectrum of $L^2(G/H_{ss})$ are automorphic representations [B-S1]. Most of these automorphic representations are known and have been constructed using other techniques  for example in [K-R],  [H-P], [S], [M-W].

 \medskip
 \subsection{Example 4}  $G=SL( 2n,\bR)$, We take $H$ as the connected component of $S(GL(p,\bR)\times GL(q,\bR)$. 
Then $H_{ss}=SL(p,\bR) \times SL(q,\bR)$  where $p+q = 2n$.  
 and  \[ \mbox{rank } G/H  = \mbox{rank } K/K\cap H=\mbox{min }(p,q).\]
 
The results of E. van den Ban and H. Schlichtkrull show that all the representations in the continuous 
spectrum are unitarily  induced from $\theta \tau$ - stable parabolic subgroups. It is easy to see that these parabolic subgroups are all cuspidal and thus  the representations in the discrete spectrum determine whether $L^2(SL(n,\bR)/H)$ is tempered. 
 
 We recall the parametrization of the representations in the discrete spectrum. 
Using the decomposition $\bh\otimes \bC \oplus \bq \otimes \bC $ of $gl(2n,\bC)$  
we conclude that the skew diagonal matrices in $\bq\otimes \bC$  are a maximal abelian 
subspace of $so(2n,\bC) \cap \bq\otimes \bC$ . By [M-O] their centralizer  $L$ is the Levi 
subgroup of a $\theta$ stable parabolic subgroup. The  representations in the discrete spectrum 
are cohomologically induced from a character of the subgroup $L$. If the commutator subgroup $L$ doesn't contain a noncompact semi simple subgroup then the representations are tempered.([K-V] Chap. XI or [Z-V] 6.16). Thus we conclude:
 
  \begin{itemize} 
\item If $p = q = n$ the subgroup $[L,L]$ is a product of $n$ compact tori. Hence all representation 
$L^2(SL(2n,\bR)/(SL(n,\bR)\times SL(n,\bR)) $ in the discrete spectrum are tempered and thus  $L^2(SL(2n,\bR/H_{ss})$  is tempered. 

\item
If $p-q \geq 2$ then $L$ has a non compact subgroup and hence   the representations in the discrete spectrum of $L^2(G/H_{ss})$ are the Langlands subquotient of representations which is not unitarily induced.
Hence $ L^2(SL(2n,\bR)/(SL(n,\bR)\times SL(n,\bR)) $ is not tempered.

\item
Using remark 2 in section III we can construct a large number of additional semi-simple subgroups $H_{ss}$ so that $L^2(SL(2n,\bR/H_{ss})$  is tempered.
\end{itemize}
The work of C. Moeglin and J.L. Waldspurger [M-W] shows that these representations 
are in the residual spectrum of a congruence subgroup of $GL(n,\bR).$ 
Similar considerations for general linear groups can be found in [V].

 \medskip
\subsection{Example 5}
$G= SO(p,q)$,  $p+q =2n \geq 4 $ with $p \geq q > 2$ and $H = SO(1,1) \times SO(p-1, q-1)$ and $H_{ss}= SO(p-1,q-1)$.  

\medskip \noindent
{\bf Claim:}
\begin{itemize}
\item
$L^2(SO(p,q)/SO(p-1,q-1))$  is not tempered.
\end{itemize}
We have 
 \[ \mbox{rank } G/H  = \mbox{rank } K/K\cap H=2\]
We argue as in example 2. Then $[L,L]$ has a factor isomorphic to $S0(p-2,q-2)$ and is hence noncompact. 
So there are non tempered representations in the discrete spectrum.

  In [K] T. Kobayashi  considered the case   $G/H_0$ where $H= H_c\times H_{0}$ where $H_c$ is compact orthogonal group and $H_{0}$ is a noncompact orthogonal group and determined the 
 multiplicities in the discrete spectrum.

\medskip
\subsection{Example 6}  $G= Sp(n,\bR )$, $H= GL(n,\bR)$, $H_{ss}=SL(n,\bR )$.  

\medskip
\noindent
{\bf Claim:}
\begin{itemize}

\item $L^2(Sp(n,\bR)/SL(n,\bR) $ is tempered.
\end{itemize}

The proof proceeds as follows: \\
\underline {1. Step:} All the representations in the discrete spectrum are tempered.\\
\underline{2. Step:}  Each conjugacy class of parabolic subgroups contains a $\theta \tau $--invariant parabolic subgroup $MAN$.\\
\underline{3. Step:} All discrete series representations of $M/M\cap H$ are tempered.

For simplicity assume that the symplectic group is defined by the form quadratic form defined by the matrix
\[
\left(
\begin{array}{cc}
0  &   I   \\
-I  &   0   \\ 
\end{array}
\right)
\]
where $I$ is the identity matrix.
The subgroup $H=GL(n,\bR )$ of $G $  is the fix point set of the automorphism $\tau$ defined by conjugation with
\[
\left(
\begin{array}{cc}
 I & 0   \\
 0 & I     
\end{array}
\right)
.\]
The  maximal  compact subgroup $K_H$ of $H$ is  $K \cap H=O(n)$.
           Furthermore 
            \[ \bg= \bk \oplus \mathfrak p\]
 \[ \bg = \bh \oplus \bq\] 
The one-dimensional torus $T_0$ in the center of  $K$ also defines a torus on $K/H\cap K$.
  Its Lie algebra $\bt_0$ is direct summand of the maximal abelian subalgebra $\ba_k$ of $\bq_k=\bk \cap \bq$. Since $T_0$ defines the complex structure on the symmetric space $G/K$ the centralizer of $\ba_k$  in $G$ is contained in $K$. Thus every representation in the discrete spectrum is tempered.
  
 The   $\theta\tau $--stable parabolic subalgebras are determined by maximal  abelian subspaces i
 in $\mathfrak p \cap \bq$. Now
 \[\tilde{\bh}:=\bk \cap \bh \oplus \mathfrak p\cap \bq = gl(n,\bR)\]
 is the fix point set of the involution $\theta\tau$ since the fix point set of $\theta\tau$ is conjugate in GL(2n, $\bR$) to GL(n,$\bR$).
 This implies that there is a n dimensional abelian split subalgebra  $\tilde{\ba}_H$ in  $p \cap q$ consisting of the matrices
 \[
\left(
\begin{array}{cc}
0  &  D   \\
 D   &  0  
\end{array}
\right)
\]
where D is real diagonal matrix.
  Hence every  conjugacy classes of parabolic subgroups contains a
$\theta \tau$ -stable parabolic $P_s=M_sA_sN_s$ whose Levi subgroup is a centralizer of  $\tilde{\ba}_H$.
  
Next we have to determine $M_s\cap H$. Note that $M_s$ is a product of general linear groups and a symplectic 
group. The factors isomorphic to general linear groups are subgroups of $\tilde{H}$. 
Since $\tilde{H} \cap H = K \cap H$ is an orthogonal group, the intersection of the general linear subgroups  of $M_s$ with $H$   are orthogonal groups and hence the corresponding symmetric space has no discrete spectrum. Thus we may assume that $M_s=Sp(m,\bR)$ with $m<n$. In this case $\theta\tau$ is an involution of $M_s$ with fix points $\tilde{H} \cap M_s= GL(m,\bR)$.  Furthermore since $\theta $ and $ \tau $ commute their restriction to $M_s$    also defines an automorphisms of $M_s$. So the fix point set of $\theta \tau _{|M_s}$ is conjugate to to the fix point set of $\tau_{|M_s}$ in $GL(2n,\bR)$. Hence 
we conclude that $M_s \cap {H} $ is isomorphic to $GL(m,\bR).$
By step 1 the representations in  the discrete spectrum of $Sp(m,\bR)/GL(m,\bR)$ are tempered and thus by [B-S] all the representations in the continuous spectrum of $Sp(n,\bR)/SL(n,\bR)$ are tempered.

\medskip
\subsection{Example 7}  Cayley type spaces considered in Olafson Orsted [O-O2] and [F-K].
These are  \\
a.) $G=Sp(n,\bR)$, $H=GL(n,R)$ and $H_{ss}=SL(n,\bR )$, $n>1$\\
b.) $G=SO(2,n)$, $H=SO(1,1)SO(1,n-1) $ and $H_{ss} = SO(1,n)$, $n>2$\\
c.) $ G=SU(n,n)$, $H=SL(n,C) \bR^+$ and $H_{ss}=SL(n,\bC)$\\
d.) $ G= O^*(2n)$,  $H=\bR^+SU^*(2n) $ and $H_{ss} = SU^*(2n)$\\
e.) $G= E_7(-25)$, $H=E_6(-26)\bR^+ $  and $H_{ss} = E_6(-26)$\\

\medskip \noindent
{\bf Claim:}

\begin{itemize}
\item
In examples 5b, c, d  and  $n$ large enough $L^2(G/H_{ss})$ is not tempered.
\item 
In example 5a $L^2(G/H_{ss})$ is tempered.
\item
We expect that in example 5e  $L^2(G/H_{ss})$ is tempered.
\end{itemize}

The proof is based on case by case considerations of  the spectrum of  $L^2(G/H)$. In  [O-O2] it is proved that all these spaces are of equal rank and hence $L^2(G/H)$ has a discrete spectrum.

Case b.)  The arguments in example 5 show that  the  representations in the discrete spectrum of $L^2(SO(n,2)/SO(n-1,1))$ are tempered iff $n\leq 2$.  So we can conclude that 
$L^2(SO(n,2)/SO(n-1,1))$ is not tempered if $3\leq n$.

Case c.)  It was proved in [O-O1] that the discrete spectrum for $SU(n,n)/H )$ contains some  non tempered highest weight representations. Hence $L^2(SU(n,n)/S L(n,\bC))$ is not tempered.

Case a.) This was proved in example 6.

Case d.) We have $rank(G/H) = n$.  The Levi of the $\theta$-- stable parabolic subgroup also contains a subgroup of type $A_{2n-1}$. Since it is not the maximal compact subgroup, $L$ has a noncompact  subgroup. This implies  that the discrete spectrum of  $L^2(O^*(2n)/SU^*(2n))$ is not tempered.

Case e.) We only prove that the discrete spectrum is tempered. The arguments are the same as in  example 4. Recall that
\begin{enumerate}
               \item the rank of $G/H$ is 3 
                \item  the maximal  compact compact subgroup $K$ of $G $ is $E_6 SO(2)$
              \item  the maximal compact subgroup $K_H$ of $H$ is  $F_4$
              \item $K/H\cap K$ has a one-dimensional  compact torus $T_0$ as factor. 
 \end{enumerate}
   The centralizer of this torus $T_0$ is $K$. Its Lie algebra is direct summand of the maximal abelian subalgebra $\ba_k$ of $\bq_k=\bk \cap \bq$. Since $T_0$ defines the complex structure on the symmetric space $G/K$ the centralizer of $\ba_k$  in $G$ is contained in $K$. Thus every representation in the discrete spectrum is tempered.
 
 As in example 6 we conclude that the fix point set of $\theta \tau $ is a subgroup 
isomorphic to  $H=E_6(-26)\bR^+ $, which has real rank 3.   Hence every  conjugacy classes of parabolic subgroups contains a
$\theta \tau$ -stable parabolic $P_s=M_sA_sN_s$ whose Levi subgroup is a centralizer of  $\tilde{\ba}_H$.

\medskip

\subsection{Example 8} $G= SL(2n, \bC)$ and $H_{ss} $ has a covering $T^1SL(p,\bC) \times SL(q,\bC)$, p+q = 2n for a one dimensional torus $T^1$. Then
\[L^2(SL(n,\bC)/SL(p,\bC)\times SL(q,\bC))= \oplus _{\delta \in  \hat{T}}L^2(SL(n,\bC)/H_{ss},\delta)\]
where $L^2(SL(n,\bC)/H_{ss},\delta)$ are the $L^2$-section of the line bundle defined by the character $\delta $ of $H_{ss}$.
As in example 2 we are in the equal rank case. 

The same arguments as in example 4 show 

\medskip \noindent
{\bf Claim:}
   \begin{itemize} 
\item If $p = q = n$ the subgroup $[L,L]$ is compact. Hence all representations  in the discrete spectrum of $ L^2(SL(2n,\bC)/H_{ss}) $ are tempered and thus  $L^2(SL(2n,\bC/H_{ss})$  is tempered  

\item
If $p-q \geq 2$ then $[L,L]$ is not compact and hence   the representations in the discrete spectrum of $L^2(G/H_{ss})$ are the Langlands subquotient of representations which is not unitarily induced .
Hence $ L^2(SL(2n,\bC)/SL(n,\bC)\times SL(n,\bC)) $ is not tempered.
\end{itemize} 
\medskip

\section{Bibliography} 
\begin{description} 
\item[B-K] \ Y. Benoist, T. Kobayshi, Temperedness of reductive homogeneous spaces, to appear in {\it European Journal of Mathematics},
\item[B-S] \ E. van den Ban, H. Schlichtkrull,  The Plancherel formula for a reductive symmetric space II: Representation theory, Invent. Math 161, 2005, pages  567--628.
\item[B-S1] \  M. Burger, P. Sarnak, Ramanujan Duals II, Inventiones
Math. 106, 1991, pages 1--11. 
\item
[F-K] \ J. Faraut, A. Koranyi,  {\it Analysis on Symmetric cones}, Oxford Mathematical Monographs, 1994 
\item [F] \ J. M. G. Fell, Weak containment and induced representations of groups. Canad. J. Math.,
14, 1960, pages 237--268.
\item
[H-P] \ 
R. Howe, I.  Piatestski-Shapiro; A counter example to the the "generalized Ramanujan conjecture", {\it Automorphic forms, Representations and L-functions, Part 1,} AMS conference series 1979.
\item [K-V] \ A.W Knapp, D. Vogan,  {\it Cohomological Induction and Unitary Representations}, Princeton University Press 1996.
\item
[K] \ T.  Kobayashi, {\it Singular unitary representations and discrete series for indefinite Stiefel manifolds U(p,q;F)/U(p-m,q;F)}, Mem. Amer. Math. Soc., vol. 462, Amer. Math. Soc., 1992.
\item [K-R]  \ S. Kudla and S. Rallis, Degenerate principal series and invariant distributions, Israel J. Math., 69, 1990, pages. 25-45.
\item[M-W] \ C. Moeglin, J.-L. Waldspurger, 
Le spectre rŽsiduel de GL(n), Ann. Sci. ƒcole Norm. Sup.  22, 1989, 605--674.
\item
[O-O1] \ G. Olafson and B. {\O}rsted, The holomorphic discrete series for affine symmetric spaces,
J. Funct. Anal. 81, 1988, pages 126--159.
\item
[O-O2] \ G. Olafson and B. {\O}rsted,  Causal compactification and Hardy spaces, Transactions of the AMS vol 351,(9), pages 3771-3792
\item [O-M] \ T. Oshima and T. Matsuki, A Complete Description of Discrete Series for Semisimple Symmetric Spaces, Adv. Stud. Math. 4, 1984,  pages  332-390.
\item [S] \ H. Schlichtkrull, The Langlands parameters of Flensted-Jensen's discrete series for semisimple symmetric spaces, Journal of Functional Analysis,  
Volume 50, Issue 2, 1983, pages 133Ð-150.
\item[Sp] \ B. Speh, Unitary representations of GL(n,$\bR$) with non-trivial  $(\bg,K)$ -cohomology,
Invent. Math. 71, 1983, pages 443-465.
\item  [V-Z] \ G. Zuckerman, D. Vogan, Unitary representations with nonzero  $(\bg,K)$ -cohomology,
Composition Math, 53, 1984, pages 51--90.
\item [V] \ A. Venkatesh, The Burger-Sarnak method and operations on the unitary dual of GL(n),
 Represent. Theory 9,  2005, pages 268--286.
\item[W] Wallach, Nolan R., {\it Holomorphic continuation of generalized
Jacquet integrals for degenerate principal series},
Represent. Theory  10  (2006), 380--398 (electronic).

\end{description}
\end{document}